\documentclass{amsart}

\usepackage{amssymb}

\theoremstyle{definition}

\theoremstyle{remark}

\numberwithin{equation}{section}
\newcommand{\R}{\mathbb{R}}  

\begin{document}

\begin{center}
{\large The impact on mathematics of the paper

\smallskip

{\it Oscillation and chaos in physiological control systems}

\smallskip

by Mackey and Glass in Science, 1977}
\end{center}

\bigskip

\bigskip

\begin{center}
Hans-Otto Walther
\end{center}

\bigskip

\begin{center}
November 2009
\end{center}

\bigskip

\bigskip






The most important examples for the study of autonomous nonlinear
delay differential equations, which describe delayed feedback,
during the last decades have been Wright's equation
$$
x'(t)=-\alpha\,x(t-1)[1+x(t)]
$$
from 1955
\cite{Wr} and the Mackey-Glass equations
\begin{equation}
P'(t)=\frac{\beta_0\theta^n}{\theta^n+(P(t-\tau))^{n}}-\gamma\,P(t)
\end{equation}
and
\begin{equation}
P'(t)=\frac{\beta_0\theta^nP(t-\tau)}{\theta^n+(P(t-\tau))^{n}}-\gamma\,P(t)
\end{equation}
with positive parameters $\beta_0,\theta,\gamma,\tau$ and an integer
$n>0$, from 1977 \cite{MG}. Closely related to Eq. (0.2) is
\begin{equation}
x'(t)=-p\,x(t)+e^{-x(t-h)}(x(t-h))^8
\end{equation}
with parameters $p>0$, $h>0$. According to my knowledge this form
was proposed by Andrzej Lasota \cite{L0}. All of these equations can
be transformed to
\begin{equation}
x'(t)=-\mu\,x(t)+f(x(t-1))
\end{equation}
with a parameter $\mu\ge0$ (for Wright's equation as long as
solutions $x>-1$ are considered). The main difference is that in
case of Wright's equation and in Eq. (0.1) the nonlinearities $f$
are monotone while for the Mackey-Glass equation (0.2) and for Eq.
(0.3) $f$ is unimodal, or a hump function. This causes entirely
different solution behaviour: In case of monotonicity there is a
planar attractor of almost all solutions, with nested periodic
orbits, see \cite{W4} and also \cite{M-PS1,M-PS2}, while for
(suitable) hump functions chaotic dynamics exist (which are
impossible for planar flows).

\medskip

In the sequel I describe the role of the Mackey-Glass equation
(0.2) and of Eq. (0.3) in the study of nonlinear delay
differential equations, as far as rigorous mathematical results
are concerned, and from the very personal point of view of my own
involvement. Numerical work will only briefly be touched. In fact,
there is a much wider area of work in more applied fields which
has been inspired by the Mackey-Glass equations.

\medskip

The latter were proposed by Michael Mackey and Leon Glass as a model
for the feedback control of blood cells, among others. At about the
same time, even a little earlier, another equation of the form (0.1)
with a monotone nonlinearity had been suggested for the same
purpose, by Lasota and Marie Wazewska-Czyzewska \cite{LWC}.

\medskip

Immediately after the paper \cite{MG} had appeared it attracted
the attention of mathematicians interested in nonlinear dynamics
and functional differential equations. Certainly the word {\it
chaos} coined by Jim Yorke in 1975 for complicated solution
behaviour played a role here. The paper \cite{LY} by Tien-Yien Li
and Yorke on interval maps had popularized the theme of
complicated motion, after the work of Steven Smale \cite{Sm} and
Leonid Shilnikov \cite{Shi} on diffeomorphisms and flows (and
after Aleksandr Sharkovsky's beautiful paper on interval maps
\cite{Sh} which became widely known much later).

\medskip

It was in the summer of 1978, during a 4 month visit to {\it
Sonderforschungsbereich 72} at the university in Bonn, that I
began to learn about chaotic dynamics, together with a group of
graduate students of Heinz-Otto Peitgen, among them Norbert
Angelstorf, Hubert Peters, Michael Pr\"ufer, and Hans-Willi
Siegberg. All of us were working on functional differential
equations. The main topic at the time were existence and global
bifurcation of periodic solutions, and the relations to fixed
point theory. The Li-Yorke route to chaos via period doubling
bifurcations (the Feigenbaum scenario \cite{F}) was clearly most
appealing to anyone in this area (of FDEs), with the chaos
observed by Mackey and Glass in simple-looking, nice equation as a
far aim. Nowadays mathematical results came quite close to this
aim - but it is not yet reached and remains a challenge which
calls for new ideas. - In a seminar in Bonn we became acquainted
with shift dynamics and Li-Yorke chaos for interval maps
\cite{LY}, the little brother of chaos for horseshoe
homeomorphisms \cite{Sm}.

\medskip

In 1979 there appeared Roger Nussbaum's paper \cite{N1} which
established coexistence of different periodic orbits for Eq. (2)
with $\mu=0$ and certain hump-shaped nonlinearities, as in the
Mackey-Glass equation (0.2). One of the first proven results (to
my knowledge) exactly for the Mackey-Glass equations is local Hopf
bifurcation of periodic solutions, due to Mario Martelli, Klaus
Schmitt and Hal Smith \cite{MSS}.

\medskip

In the small group in Bonn Hubert Peters rediscovered that step
functions as nonlinearities can make dynamics generated by a delay
differential equation accessible, and that one can compute
periodic solutions explicitly. (Similar ideas had been used
earlier, also in work on complicated dynamics, see e. g.
\cite{Le}, but we were not yet aware of this.) My contribution was
the observation that in case of such simple nonlinearities one can
introduce local coordinates and reduce the semiflow as well as
Poincar\'e return maps to explicitly given maps in finite
dimensions. (Much later this became chapter XVI in the monograph
\cite{DvGVLW}). In the following years Peters and Siegberg reduced
Poincar\'{e} return maps for certain delay differential equations
with step functions as nonlinearities to interval maps and
obtained Li-Yorke chaos in this way \cite{P,Si}. I showed that for
equations with smooth nonlinearities (which were still constant on
long intervals) there are unstable periodic orbits with
homoclinics, and that close to this structure Li-Yorke chaos
exists \cite{W1}.

\medskip

During a year at Michigan State University, 1979-1980, I met Lasota
and Pavol Brunovsky - all of us being brought together by Shui-Nee
Chow. From Andy I learnt about his work with Wazewska-Czyzewska
\cite{LWC}, and about his topological and measure-theoretic ideas on
chaos which gave us another strong impetus to explore complicated
motion. Also his interest in applications in physiology and cell
biology was very inspiring. Not being experts in stochastics, Andy's
paper \cite{L} was closest to our abilities. Palo got into the area
and obtained a result for first order PDEs \cite{B} which was
related to \cite{L}. - A little later Uwe an der Heiden and I
obtained Li-Yorke chaos for a delay differential equation with a
smooth nonlinearity (but still close to a step function) which was
modelled after the Mackey-Glass equation (0.2) and Eq. (0.3)
\cite{adHW}.

\medskip

Another area of research which to some extent was motivated by the
route to chaos in the Mackey-Glass equations concerns the
development of numerical methods for bifurcation diagrams which
follow the fixed point index along the branches in the diagram.
J\"urgens, Peitgen, Pr\"ufer and Dietmar Saupe \cite{Sa} worked in
this direction, both on the theoretical foundation and on the
improvement of algorithms.

\medskip

I became interested in the very first step of this route and
obtained a result on secondary bifurcation of periodic orbits in
delay differential equations with hump nonlinearities, like the
{\it sine} function \cite{W2}. This included a new approach to
Floquet multipliers of periodic solutions of such equations, via a
characteristic equation which involved a boundary value problem
for an ODE. Since then this approach to Floquet multipliers has
been used and extended in numerous papers, and is still the object
of on going research. Work with Chow \cite{CW} on stability and
more recently with Alexander Skubachevsky \cite{SkW1,SkW2}, and by
Nikolay Zhuravlev \cite{Z,ZSk} must be mentioned here.

\medskip

Nussbaum \cite{N2}and Steven Chapin \cite{Ch1,Ch2} extended the
nonuniqueness result from \cite{N1} and achieved very sharp
estimates of the shape of periodic solutions of Eq. (0.2) with
$\mu=0$.

\medskip

In a series of papers \cite{M-PN1,M-PN2,M-PN3,M-PN4,M-PN5} John
Mallet-Paret and Nussbaum found that the routes to chaos for the
parametrized difference equation
\begin {equation}
x_{n+1}=f(x_n)
\end{equation}
and for the delay differential equation
\begin{equation}
\epsilon\,\dot{x}(t)=-x(t)+f(x(t-1))
\end{equation}
with $\epsilon>0$ small are not the same; the dynamics given by
Eq. (0.6) is much richer, and part of it is mirrored by a
difference equation which is different from Eq. (0.5).

\medskip

Relations between Eq. (0.4) and the difference equation
$$
x(t)=f(x(t-1)),\,\,t\in\R,
$$
were also discussed in the survey paper \cite{IS} by Anatoli Ivanov
and Sharkovsky.

\medskip

Mallet-Paret \cite{M-P} constructed a Morse decomposition of the
global attractor for equations
\begin{equation}
x'(t)=f(x(t),x(t-1))
\end{equation}
which generalize Eq. (0.4) with negative and positive feedback. The
invariant Morse sets of this decomposition may contain chaotic
solutions.

\medskip

A challenge in the early 1980ies certainly was to get existence of
chaos for delay differential equations with nonlinearities not
resembling step functions, and not only Li-Yorke chaos which is tied
to non-injectivity but the kind of chaos for injective horseshoe
maps where a Bernoulli shift in finitely many symbols represents the
complicated solution behaviour. A first result in this direction is
\cite{W3}, for an equation without decay term ($\mu=0$ in Eq. (0.4)
above) and for nonlinearities which are periodic (not in time, of
course). The structure underlying chaos is a chain of unstable
periodic orbits with heteroclinic connections. The result in
\cite{W3} required the extension of general theory on hyperbolic
sets, shadowing and transversal homoclinic points from the case of
diffeomorphisms in finite dimensions to the case of continuously
differentiable maps, not necessarily injective, in arbitrary Banach
spaces. This extension was done in work with Heinrich Steinlein
\cite{SW1,SW2}, who proposed the crucial wider notion of hyperbolic
sets which made this theory elegant. Further extensions of the
general theory are due to Bernhard Lani-Wayda \cite{L-W1}, who also
obtained existence of chaotic motion for delay differential
equations in a Shilnikov scenario, with a homoclinic loop at a
stationary state and with heteroclinics connecting stationary states
\cite{L-W3}.

\medskip

Not long ago Lani-Wayda and Roman Szrednicki established existence
of chaotic motion for delay differential equations by means of a
generalized Lefschetz fixed point theorem \cite{L-WS}.

\medskip

Considering the route to chaos via bifurcations of periodic
solutions it must be said that period doubling could not yet be
shown in case of delay differential equations. However there are
deep and beautiful results on sequences of secondary bifurcations
along a branch of periodic solutions, for equations (0.4) without
decay term. This is work of Peter Dormayer \cite{D1,D3,D4}. Related
is work by him on Floquet multipliers \cite{D2}, and jointly with
Lani-Wayda \cite{DL-W}, also on Floquet multipliers and secondary
bifurcation of periodic solutions.

\medskip

Christopher McCord and Konstantin Mischaikow constructed a
semiconjugacy between the flow on the attractor of a delay
differential equation and a model flow of an ordinary differential
equation, using ideas related to the Conley index \cite{MM}.

\medskip

A further step towards a proof of chaotic motion in the Mackey-Glass
equation was done in the 1990ies. Lani-Wayda and I obtained chaos
for Eq. (0.3) with $\mu=0$ in the negative feedback case
\cite{L-WW1,L-WW2}. Here chaos occurs near a solution which is
homoclinic to a hyperbolic unstable periodic orbit. The latter was
found in joint work with Ivanov \cite{IL-WW}. The nonlinearities
considered in \cite{L-WW1,L-WW2} are smooth and not close to step
functions but have a complicated shape with several humps. Later
Lani-Wayda improved the technique and succeeded to establish chaos
for equations with nicer nonlinearities, with only one hump - as in
the Mackey-Glass equation \cite{L-W2}. Up to now this is the result
which is closest to a proof of chaos in the Mackey-Glass equations.

\medskip

Conditions for stability and attractivity of the positive stationary
state of Eqs. (0.2) and (0.3) were studied by Ivanov, Eduardo Liz,
Clotilde Martinez, Victor Tkachenko, Elena and Sergei Trofimchuk,
and Gergely R\"ost \cite{ILT,LMT,LTkTr,LTT,LR,R}.

\medskip

Closely related to the Mackey-Glass equation is work by Mackey,
Chunhua Ou, Laurent Pujo-Menjouet, and Jianhong Wu \cite{MOP-MW} on
existence and stability of periodic solutions with long period, for
a delay differential system which is a more specific model for the
regulation of blood cell populations.


\bigskip

\begin{thebibliography}{999}

\bibitem{adHW} an der Heiden, U., and H. O. Walther, {\it Existence of chaos in control systems
with delayed feedback.} Journal of Differential Equations 47
(1983), 273-295.


\bibitem{B} Brunovsky, P., {\it Notes on chaos in the cell population partial differential equation.}
Nonlinear Analysis 7 (1983), 167-176.


\bibitem{Ch1} Chapin, S. A., {\it Asymptotic analysis of
differential-delay equations and nonuniqueness of periodic
solutions.} Math. Methods in the Applied Sciences 7 (1985), 223-237.


\bibitem{Ch2} Chapin, S. A., {\it Periodic solutions of a
differential-delay equation and the fixed point index.} Applicable
Analysis 33 (1989), 119-126.


\bibitem{CW} Chow, S. N., and H. O. Walther, {\it Characteristic multipliers and stability of
symmetric periodic solutions of $\dot{x}(t)=g(x(t-1)))$.}
Transactions of the A. M. S. 307 (1988), 127-142


\bibitem{DvGVLW} Diekmann, O., van Gils, S. A., Verduyn Lunel, S. M., and H. O. Walther, {\it
Delay Equations: Functional-, Complex- and Nonlinear Analysis}.
Springer, New York, 1995.


\bibitem{D1} Dormayer, P., {\it Smooth bifurcation of symmetric periodic solutions of functional differential equations.}
J. Differential Equations  82 (1989), 109-155.


\bibitem{D2} Dormayer, P., {\it An attractivity region for characteristic multipliers of special
symmetric solutions of $\dot{x}(t)=\alpha\,f(x(t-1))$.}
J. Math. Analysis and Applications 168 (1992), 70-91.


\bibitem{D3} Dormayer, P., {\it Smooth symmetry-breaking bifurcation for functional differential equations.}
Differential and Integral Equations 5 (1992), 831-854.


\bibitem{D4} Dormayer, P., {\it Floquet multipliers and secondary bifurcation of periodic solutions
of functional differential equations.}
Habilitation thesis, Gie{\ss}en, 1996.


\bibitem{DL-W} Dormayer, P., and B. Lani-Wayda, {\it Floquet
multipliers and secondary bifurcations in functional differential
equations. Numerical and analytical results.} J. Applied Math. and
Physics 46 (1995), 823-858.


\bibitem{F} Feigenbaum, M. J., {\it The universal metric properties of nonlinear transformations.}
Journal of Statistical Physics 19 (1978), 25-52.


\bibitem{HMPMW} Hu, C. O., Mackey, M. C., Pujo-Menjouet, L., and
J. Wu, {\it Periodic oscillations of blood cell populations in
chronic myelogeneous leukemia.} SIAM J. Math. Analysis 38 (2006),
166-187.


\bibitem{IL-WW} Ivanov, A. F., Lani-Wayda, B., and H. O. Walther, {\it Unstable hyperbolic
periodic solutions of differential delay equations.} In {\it
Recent Trends in Differential Equations}, R.P. Agarwal ed.,
301-316, WSSIAA vol 1, World Scientific, Singapore, 1992.


\bibitem{ILT} Ivanov, A. F., Liz, E., and S. Trofimchuk, {\it Global
stability of a class of nonlinear delay differential equations.}
Differential Equations and Dynamical Systems 11 (2003), 33-54.


\bibitem{IS} Ivanov, A. F., and A. N. Sharkovsky, {\it Oscillations in singularly perturbed delay equations.}
In {\it Dynamics Reported}, Jones, C.K.R.T., Kirchgraber, U., and H.
O. Walther eds., Vol. 1 (new series), pp. 164-224, Springer, New
York, 1992.


\bibitem{L-W1} Lani-Wayda, B., {\it Hyperbolic sets, shadowing and persistence for nonlinear mappings in Banach spaces.}
Pitman Research Notes in Math., No. 334 (1995).


\bibitem{L-W2} Lani-Wayda, B., {\it Erratic solutions to simple delay equations.} Trans. A. M. S. 351 (1999), 901-945.


\bibitem{L-W3} Lani-Wayda, B., {\it Wandering solutions of delay
equations with sine-like feedback.} Memoirs of the A. M. S., No.
718 (2001).


\bibitem{L-WS} Lani-Wayda, B., and R. Szrednicki, {\it A generalized Lefschetz fixed point theorem and symbolic dynamics
in delay equations.} Ergodic Theory and Dynamical Systems 22 (2002),
1215-1232.


\bibitem{L-WW1} Lani-Wayda, B., and H.O. Walther, {\it Chaotic motion generated by delayed
negative feedback. Part I: A transversality criterion.}
Differential and Integral Equations 8 (1995), 1407-1452.


\bibitem{L-WW2} Lani-Wayda, B., and H. O. Walther, {\it Chaotic motion generated by delayed
negative feedback. Part II: Construction of nonlinearities.}
Mathematische Nachrichten 180 (1996), 141-211.


\bibitem{L0} Lasota, A., personal communication, 1980.


\bibitem{L} Lasota, A., {\it Stable and chaotic solutions of a
first-order partial differential equation.} Nonlinear Analysis 5
(1981), 1181-1193.


\bibitem{LWC} Lasota, A., and M. Wazewska-Czyzewska, {\it
Matematyczne problemy dynamiki ukladu krwinek czerwonych.} Mat.
Stosowana 6 (1976), 23-40.


\bibitem{Le} Levinson, N., {\it A second order differential equation
with singular solutions.} Annals of Math. 50 (1949), 126-153.


\bibitem{LY} Li, T. Y., and J. Yorke, {\it Period three implies
chaos}, American Math. Monthly 82 (1975), 985-992.


\bibitem{LMT} Liz, E., Martinez, C., and S. Trofimchuk, {\it
Attractivity properties of infinite-delay Mackey-Glass type
equations.} Differential and Integral equations 15 (2002), 875-896.


\bibitem{LR} Liz, E., and G. R\"ost, {\it On global atractors for
delay differential equations with unimodal feedback.} Discrete and
Continuous Dynamical Systems, to appear.


\bibitem{LTkTr} Liz, E., Tkachenko, V., and S. Trofimchuk, {\it A
global stability criterion for scalar functional differential
equations.} SIAM J. Math. Analysis 35 (2003), 596-622.


\bibitem{LTT} Liz, E., Trofimchuk, E., and S. Trofimchuk, {\it
Mackey-Glass type delay differential equations near the boundary
of absolute stability.} J. Math. Analysis and Applications 275
(2002), 747-760.


\bibitem{MG} Mackey, M. C., and L. Glass, {\it Oscillation and chaos in physiological control systems.}
Science 197 (1977), 287-295.


\bibitem{MOP-MW} Mackey, M. C., Ou, C., Pujo-Menjouet, L.,  and J.
Wu, {\it Periodic oscillations of blood cell populations in chronic
myelogenous laukemia.} SIAM J. Math. Analysis 38 (2006), 166-187.


\bibitem{M-P} Mallet-Paret, J., {\it Morse decompositions for
delay-differential equations.} J. Differential Equations 72 (1988),
270-315.


\bibitem{M-PN1} Mallet-Paret, J., and R. D. Nussbaum, {\it Global
continuation and complicated trajectories of periodic solutions
for a differential-delay equation.} In {\it Proc. of Symposia in
Pure Math.}, Vol. 45, Part 2, pp.155-167, A. M. S., Providence
(RI), 1986.


\bibitem{M-PN2} Mallet-Paret, J., and R. D. Nussbaum, {\it Global
continuation and asymptotic behaviour for periodic solutions of a
differential-delay equation.} Annali di Matematica Pura ed Applicata
145 (1986), 33-128.


\bibitem{M-PN3} Mallet-Paret, J., and R. D. Nussbaum, {\it A
bifurcation gap for a singularly perturbed delay equation.} In {\it
Chaotic Dynamics and Fractals}, Barnsley, M., and S. Demko eds., pp.
263-287, Academic Press, New York, 1986.


\bibitem{M-PN4} Mallet-Paret, J., and R. D. Nussbaum, {\it A
differential-delay equation arising in optics and physiology.}
SIAM J. Math. Analysis 20 (1989), 249-292.


\bibitem{M-PN5} Mallet-Paret, J., and R. D. Nussbaum, {\it Multiple transition
layers in a singularly perturbed differential-delay equation}.
Proceedings of the Royal Society of Edinburgh 123 A (1993),
1119-1134.


\bibitem{M-PS1} Mallet-Paret, J., and G. Sell, {\it Systems of delay
differential equations: Floquet multipliers and discrete Lyapunov
functions.} J. Differential Equations 125 (1996), 385-440.


\bibitem{M-PS2} Mallet-Paret, J., and G. Sell, {\it The Poincar\'e-Bendixson theorem for
monotone cyclic feedback systems with delay.}  J. Differential
Equations 125 (1996), 441-489.


\bibitem{MSS} Martelli, M., Schmitt, K., and H. Smith, {\it
Periodic solutions of some nonlinear delay differential
equations.} J. Math. Analysis and Applications 74 (1980), 494-503.


\bibitem{MM} McCord, C. M., and K. Mischaikow, {\it On the global
dynamics of attractors for scalar delay differential equations.}
J. A. M. S. 9 (1996), 1095-1133.


\bibitem{N1} Nussbaum, R. D., {\it Uniqueness and non-uniqueness
for periodic solutions of $\dot{x}(t)=-g(x(t-1))$.}, J.
Differential Equations 34 (1979), 25-54.


\bibitem{N2} Nussbaum, R. D., {\it Asymptotic analysis of functional
differential equations and solutions of long period.} Archive for
Rational Mechanics and Analysis 81 (1983), 373-397.


\bibitem{P} Peters, H., {\it Chaotic behavior of nonlinear differential-delay equations.}
Nonlinear Analysis 7 (1983), 1315-1334.


\bibitem{R} R\"ost, G., {\it On the global attractivity controversy for a delay model of hematopoiesis.}
Applied Math. and Computation 190 (2007), 846-850.


\bibitem{Sa} Saupe, D., {\it Global bifurcation of periodic solutions to some autonomous differential delay equations.}
Applied Math. and Computation 13 (1983), 185-211.


\bibitem{Sh} Sharkovsky, A. N., {\it Coexistence of cycles of a continuous map of a line into itself.}
Ukranichkii Matematichevskii Zhurnal 16 (1964), 61-71.


\bibitem{Shi} Shilnikov, L. P., {\it On a Poincar\'e-Birkhoff problem.} Math. U. S. S. R. Sbornik 3 (1967),
353-371.


\bibitem{Si} Siegberg, H.-W., {\it Chaotic behavior of a class of differential-delay equations.}
Annali di Matematica Pura ed Applicata 138 (1984), 15-33.


\bibitem{SkW1} Skubachevsky, A. L., and H. O. Walther, {\it On Floquet multipliers of slowly
oscillating periodic solutions of nonlinear functional
differential equations.} Trudy Moskov. Mat. Obshch. 64 (2002),
3-54; English translation in: Transactions of the Moscow
Mathematical Society 2003, 1-44.


\bibitem{SkW2} Skubachevsky, A. L., and H. O. Walther, {\it On the Floquet multipliers of periodic
solutions to nonlinear functional differential equations.} Journal
of Dynamics and Differential Equations 18 (2006), 257-355.


\bibitem{Sm} Smale, S., {\it Differentiable dynamical systems}, Bulletin of the A. M. S. 73 (1967), 747-817.


\bibitem{SW1} Steinlein, H., and H.O. Walther, {\it Hyperbolic sets and shadowing for
noninvertible maps.} In {\it Advanced Topics in the Theory of
Dynamical Systems}, G. Fusco, M. Iannelli, and L. Salvadori eds.,
219-234, Academic Press, New York 1989.


\bibitem{SW2} Steinlein, H., and H. O. Walther, {\it Hyperbolic sets, transversal homoclinic
trajectories, and symbolic dynamics for $C^1$-maps in Banach
spaces.} Journal of Dynamics and Differential Equations 2 (1990),
325-365.


\bibitem{W1} H. O. Walther, {\it Homoclinic solution and chaos in $\dot{x}(t)=f(x(t-1))$}.
Nonlinear Analysis 5 (1981), 775-788.


\bibitem{W2} H. O. Walther, {\it Bifurcation from periodic
solutions in functional differential equations.} Mathematische
Zeitschrift 182 (1983), 269-289.


\bibitem{W3} H. O. Walther, {\it Hyperbolic periodic solutions, heteroclinic connections and
transversal homoclinic points in autonomous differential delay
equations.} Memoirs of the A. M. S., No. 402 (1989).


\bibitem{W4} H. O. Walther, {\it The two-dimensional attractor of $x'(t)=-\mu
x(t)+f(x(t-1))$.} Memoirs of the A. M. S. 544 (1995).


\bibitem{Wr} Wright, E. M., {\it A non-linear differential-difference equation.}
J. reine und angewandte Math. 194 (1955), 66-87.


\bibitem{Z} Zhuravlev, N. B., {\it On the spectrum of the monodromy operator for slowly oscillating
periodic solutions of functional diferential equations with several
delays.} Functional Differential Equations 13 (2006), 323-344.


\bibitem{ZSk} Zhuravlev, N. B., and A. L. Skubachevsky, {\it On
hyperbolicity of periodic solutions of functional differential
equations with several delays.} Trudy Mat. Inst. im. V. A. Steklova
256 (2007), 148-171. English translation in: Proc. Steklov Institute
of Math. 256 (2007), 136-159.


\end{thebibliography}
\end{document}